\documentclass[12pt]{amsart}
\usepackage[T1]{fontenc}
\usepackage[english]{babel}
\usepackage{graphicx}
\usepackage{amsmath,amsfonts}
\usepackage{amssymb}
\usepackage{amsthm}
\usepackage[vlined,ruled,algo2e]{algorithm2e}
\usepackage{tikz}
\usepackage{hyperref}
\usetikzlibrary{arrows,plotmarks,shapes,snakes}
\usepackage{color}

\theoremstyle{plain}
\newtheorem{theorem}{Theorem}
\newtheorem{prop}[theorem]{Proposition}

\theoremstyle{rem}
\newtheorem{exmp}{Example}

\theoremstyle{definition}
\newtheorem{defn}{Definition}

\newcommand{\setpr}[1]{{\mathbb{P}_{#1}}}

\def\P{{\mathcal{P}}}

\definecolor{light-gray}{gray}{0.65}

\begin{document}



\title[Unifying Uncertainty Representations: Generalized p-boxes]{Unifying practical uncertainty representations: I. Generalized p-boxes}


\author[S. Destercke]{S\'ebastien Destercke}
\address{Institut de Radioprotection et S\^uret\'e nucl\'eaire, B\^at 720, 13115 St-Paul lez Durance, FRANCE}
\email{sdestercke@gmail.com}
\author[D. Dubois]{Didier Dubois}
\address{Universit\'e Paul Sabatier, IRIT, 118 Route de Narbonne, 31062 Toulouse} \email{dubois@irit.fr}
\author[E. Chojnacki]{Eric Chojnacki}
\address{Institut de Radioprotection et S\^uret\'e nucl\'eaire, B\^at 720, 13115 St-Paul lez Durance, FRANCE} \email{eric.chojnacki@irsn.fr}

\begin{abstract}
There exist several simple representations of uncertainty that are
easier to handle than more general ones. Among them are random sets,
possibility distributions, probability intervals, and more recently
Ferson's p-boxes and Neumaier's clouds. Both for theoretical and
practical considerations, it is very useful to know whether one
representation is equivalent to or can be approximated by other
ones. In this paper, we define a generalized form of usual p-boxes.
These generalized p-boxes have interesting connections with other
previously known representations. In particular, we show that they
are equivalent to pairs of possibility distributions, and that they
are special kinds of random sets. They are also the missing link
between p-boxes and clouds, which are the topic of the second part of this
study.
\end{abstract}

\maketitle

\section{Introduction}

Different formal frameworks have been proposed to reason under
uncertainty. The best known and oldest one is the probabilistic
framework, where uncertainty is modeled by classical probability
distributions. Although this framework is of major
importance in the treatment of uncertainty due to variability, many
arguments converge to the fact that a single probability
distribution cannot adequately account for incomplete or imprecise
information. Alternative theories and frameworks have been proposed
to this end. The three main such frameworks, are, in decreasing
order of generality, Imprecise probability theory\cite{Walley91},
Random disjunctive sets~\cite{Dempster67, Shafer76, Molchanov05} and
Possibility theory~\cite{Zadeh78, DuboisPrade88}. Within each of
these frameworks, different representations and methods  have been proposed to make
inferences and decisions.

This study focuses on uncertainty
representations, regarding the relations existing between them, their respective expressiveness and practicality. In the past years,
several representation tools have been proposed:
capacities~\cite{Choquet54}, credal sets~\cite{levi1980a}, random
sets~\cite{Molchanov05}, possibility distributions~\cite{Zadeh78},
probability intervals~\cite{CamposAll94}, p-boxes~\cite{FersonAll03}
and, more recently, clouds~\cite{Neumaier04,NeumaierWeb04}. Such a
diversity of representations motivates the study of their
respective expressive power.

The more general representations, such as credal sets and capacities,
are expressive enough to embed other ones as particular
instances, facilitating their comparison.
However, they are generally difficult to handle, computationally
demanding and not fitted to all uncertainty theories. As for simpler
representations, they are useful in practical uncertainty analysis problems~\cite{WilliamsonDowns90,ReganAll04,FuchsNeum08}. They come in handy when trading
expressiveness (possibly losing some information) against
computational efficiency; they are instrumental
in elicitation tasks, since requesting less information\cite{BaudritDubois06};
they are also instrumental in summarizing 
complex results of some uncertainty propagation methods~\cite{FersonGinzburg95,BaudritAll06}.

The object of this study
 is twofold: first, it provides a short review of existing
uncertainty representations and of their relationships; second, it
studies the ill-explored relationships between more recent simple representations
and older ones, introducing a generalized form of p-box. Credal sets and random sets are used as the common umbrella clarifying the relations between simplified models. Finding such formal links facilitates a
unified handling and treatment of uncertainty, and suggests how
tools used for one theory can eventually be useful in the setting of
other theories. We thus regard such a study as an important and
necessary preamble to other studies devoted to
computational and interpretability issues. Such issues, which still remain a
matter of lively debate, are not the main topic of the present work,
but we nevertheless provide some comments regarding them. In
particular, we feel that is important to recall that a given 
representation can be interpreted and processed differently according to different
theories, which were often independently motivated by specific problems.

This work is made of two companion papers, one devoted to p-boxes, introducing a generalization thereof that subsumes possibility distributions. The second part considers Neumaier's clouds, an even more general representation tool.  

This paper first reviews
older representation tools, already considered by many
authors. A good complement to this first part, although adopting a subjectivist point of view, is provided by
Walley~\cite{Walley96}. Then, in Section \ref{sect:GenPboxes}, we propose and study a generalized form of p-box
extending, among other things, some results by Baudrit and
Dubois~\cite{BaudritDubois06}. As we shall see, this new
representation, which consists of two comonotonic distributions, is the missing link
between usual p-boxes, clouds and possibility distributions,
allowing to relate these three representations. Moreover, generalized p-boxes have
interesting properties and are promising uncertainty representations by themselves. In particular, Section~\ref{sec:GenPboxRS} shows that
generalized p-boxes can be interpreted as a special case of random
sets; Section~\ref{sec:GenPboxProbInt} studies the relation between
probability intervals and generalized p-boxes and discusses
transformation methods to extract probability intervals from
p-boxes, and vice-versa.

In the present paper, we restrict ourselves to uncertainty
representations defined on finite spaces (encompassing the
discretized real line) unless stated otherwise.  
Representations
defined on the continuous real line are considered in the second part of this
paper. To make the paper easier to read, longer proofs have been
moved to an appendix.


\section{Non-additive uncertainty theories and some representation tools}
\label{sec:BasicTh}

To represent uncertainty, Bayesian subjectivists advocate the use of
single probability distributions in all circumstances. However, when
the available information lacks precision or is incomplete, claiming that a
unique probability distribution can  faithfully represent
uncertainty is debatable\footnote{For instance, the following
statement about a coin: "We are not sure that this coin is fair, so
the probability for this coin to land on Heads (or Tails) lies
between 1/4 and 3/4" cannot be faithfully modeled by a single
probability.}. It generally forces to engage in too strong a
commitment, considering what is actually known.

Roughly speaking, alternative theories recalled here (imprecise
probabilities, random sets, and possibility theory) have the
potential to lay bare  the existing imprecision or incompleteness in
the information. They evaluate uncertainty on a particular event by
means of a pair of (conjugate) lower and upper measures rather than
by a single one. The difference between upper and lower measures
then reflects the lack of precision in our knowledge.

In this section, we first recall basic mathematical notions used in
the sequel, concerning capacities and the M\"obius transform.  Each
theory mentioned above is then briefly introduced, with focus on
practical representation tools available as of to-date, like
possibility distributions, p-boxes and probability intervals, their
expressive power and complexity.

\subsection{Basic mathematical notions}
\label{sec:mathtools}

Consider a finite space $X$ containing $n$ elements. Measures of uncertainty are often represented by set-functions 
called {\em capacities}, that were first introduced in Choquet's work~\cite{Choquet54}.
\begin{defn}
\label{def:capa} A capacity on $X$ is a function $\mu$, defined on
the set of subsets of $X$, such that:
\begin{itemize}
\item $\mu(\emptyset)=0,\mu(X)=1$,
\item $A \subseteq B \Rightarrow \mu(A) \leq \mu(B)$.
\end{itemize}
A capacity such that  $$\forall A,B \subseteq X, A\cap B=\emptyset,
\mu(A \cup B) \geq \mu(A) + \mu(B)$$ is said to be {\em
super-additive}. The dual notion, called {\em sub-additivity}, is
obtained by reversing the inequality. A capacity that is both {\em
sub-additive} and {\em super-additive} is called {\em additive}.
\end{defn}

Given a capacity $\mu$, its \emph{conjugate} capacity $\mu^c$ is
defined as $\mu^c(E)=\mu(X)-\mu(E^c)=1-\mu(E^c)$, for any subset $E$
of $X$, $E^c$  being its complement. In the following, 
$\setpr{X}$ denotes the set of all additive capacities on space $X$.
We will also denote $P$ such capacities, since they are equivalent
to probability measures on $X$. An additive capacity $P$ is
self-conjugate, and $P = P^c$. An additive capacity can also be
expressed by restricting it to its distribution $p$ defined on
elements of $X$ such that for all $x \in X$, $p(x) = P(\{x\})$. Then
$\forall x \in X, p(x) \geq 0$, $\sum_{x \in X} p(x) =1$ and $P(A)
=\sum_{x \in A} p(x)$.

When representing uncertainty, the capacity of a subset evaluates
the degree of confidence in the corresponding event. Super-additive and sub-additive
 capacities are fitted
to the representation of uncertainty.  The former, being
sub-additive, verify $\forall E \subset X, \mu(E) + \mu(E^c) \leq 1$ and can be called {\em cautious} capacities (since, as a consequence, 
$\mu(E) \leq \mu^c(E), \forall E$); they are tailored for modeling
the idea of certainty. The latter being
sub-additive, verify $\forall E \subset X, \mu(E) + \mu(E^c) \geq 1$, can be called {\em bold} capacities; they account for the weaker notion of
plausibility.
 
The core of a cautious capacity $\mu$ is the (convex) set of
additive capacities that dominate $\mu$, that is, $\mathcal{P}_{\mu}
= \{P \in \setpr{X} | \forall A \subseteq X, P(A) \geq \mu(A)\}$.
This set may be empty even if the capacity is cautious. We need
stronger properties to ensure a non-empty core. Necessary and sufficient conditions for non-emptiness are provided by Walley~\cite[Ch.2]{Walley91}.
However, checking that these conditions hold can be difficult in
general. 
An alternative to checking the non-emptiness of the core is to use specific
characteristics of capacities that ensure it, such as $n$-monotonicity.
\begin{defn}
A super-additive capacity $\mu$ defined on $X$ is $n-monotone$,
where $n>0$ and $n \in \mathbb{N} $, if and only if for any set
$\mathcal{A}=\{A_i | 0 < i \leq n \; A_i \subset X \}$ of events
$A_i$, it holds that
$$ \mu(\bigcup_{A_i \in \mathcal{A}} A_i) \geq \sum_{I \subseteq \mathcal{A}} (-1)^{|I|+1} \mu(\bigcap_{A_i \in I}
A_i)$$
\end{defn}
An $n$-monotone capacity is also called a Choquet capacity of order
$n$. Dual capacities are called $n$-alternating capacities. If a
capacity is $n$-monotone, then it is also ($n-1$)-monotone, but not
necessarily ($n+1$)-monotone. An $\infty$-monotone capacity is a
capacity that is $n$-monotone for every $n > 0$. On a finite space,
a capacity is $\infty$-monotone if it is $n$-monotone with $n=|X|$.
The two particular cases of $2$-monotone (also called convex)
capacities and $\infty$-monotone capacities have deserved special
attention in the
literature~\cite{ChateauneufJaffray91,Walley91,MirandaAll03}.
Indeed, $2$-monotone capacities have a non-empty core.
$\infty$-monotone capacities have interesting mathematical
properties that greatly increase computational efficiency. As we
will see, many of the representations studied in this paper possess
such properties. Extensions of the notion of capacity and of
$n$-monotonicity have been studied by de Cooman \emph{et
al.}~\cite{GDecoo05}.

Given a capacity $\mu$ on $X$, one can obtain multiple equivalent
representations by applying various (bijective)
transformations~\cite{GrabischAll00} to it. One such transformation,
useful in this paper, is the M\"obius inverse:

\begin{defn}
Given a capacity $\mu$ on $X$, its M\"obius transform is a mapping $m:
2^{|X|} \to \mathbb{R}$ from the power set of $X$ to the real line,
which associates  to any subset $E$ of $X$ the value
$$m(E) = \sum_{ B \subset E  } (-1)^{|E-B|}\mu(B) $$
\end{defn}

Since $\mu(X) = 1$,
$\sum_{E\in X}m(E)=1$ as well, and $m(\emptyset) = 0$. Moreover, it
can be shown~\cite{Shafer76} that the values $m(E)$ are non-negative
for all subsets $E$ of $X$ (hence $\forall E \in X, 1\geq m(E)\geq
0$) if and only if the capacity $\mu$ is $\infty$-monotone. Then $m$
is called a mass assignment. Otherwise, there are some
(non-singleton) events $E$ for which  $m(E)$ is negative.
Such a set-function $m$ is actually the unique solution to the set
of $2^n$ equations $$\forall A \subseteq X, \mu(A)= \sum_{E
\subseteq A} m(E),$$ given any capacity $\mu$. The M\"obius transform
of an additive capacity $P$ coincides with its distribution $p$,
assigning positive masses to  singletons only. In the sequel, we
focus on pairs of conjugate cautious and bold
capacities. Clearly only one of the two is needed
to characterize an uncertainty representation (by convention, the cautious
one).

%


\subsection{Imprecise probability theory}
\label{sec:ImpProb}

The theory of imprecise probabilities has been systematized and
popularized by Walley's book~\cite{Walley91}. In this theory,
uncertainty is modeled by lower bounds (called coherent lower
previsions)
 on the expected value that can be reached by  bounded real-valued functions on $X$ (called gambles). Mathematically speaking, such lower bounds have an expressive power equivalent to closed convex sets $\P$
of (finitely additive) probability measures $P$ on $X$.
In the rest of the paper, such convex sets will be named {\em credal
sets} (as is often done~\cite{levi1980a}). It is important to stress
that, even if they share similarities (notably the modeling of
uncertainty by sets of probabilities), Walley's behavioral
interpretation of imprecise probabilities is different from the one
of classical robust statistics\footnote{Roughly speaking, in
Walley's approach, the primitive notions are lower and upper
previsions or sets of so-called desirable gambles describing
epistemic uncertainty, and the fact that there always exists a
"true" precise probability distribution is not
assumed.}~\cite{Huber81}.

Imprecise probability theory is very general, and, from a purely
mathematical and static point of view, it encompasses all
representations considered here. Thus, in all approaches presented
here, the corresponding credal set can be generated, making the comparison of representations easier.
 To clarify this comparison, we adopt
the following terminology:
\begin{defn}
\label{def:compa} Let $\mathbb{F}_1$ and $\mathbb{F}_2$ denote two
uncertainty representation frameworks, $a$ and $b$ particular representatives of such
frameworks, and $\P_a,\P_b$ the credal sets induced by these representatives
$a$ and $b$. Then:
\begin{itemize}
\item Framework $\mathbb{F}_1$ is said to \emph{generalize} framework $\mathbb{F}_2$ if and only if
for all $b \in \mathbb{F}_2$,  $\exists a \in \mathbb{F}_1$ such that
$\P_a=\P_b$ (we also say that $\mathbb{F}_2$ is a special case of $\mathbb{F}_1$).
\item Frameworks $\mathbb{F}_1$ and $\mathbb{F}_2$ are said to be \emph{equivalent} if and only if 
for all $b \in \mathbb{F}_2$,  $\exists a \in \mathbb{F}_1$ such that
$\P_a=\P_b$ and conversely.
\item Framework $\mathbb{F}_2$ is said to be \emph{representable} in terms of framework $\mathbb{F}_1$ if and only if for all $b \in \mathbb{F}_2$, there exists a subset
$\{a_1,\ldots,a_k | a_i \in \mathbb{F}_1\}$ such that $\P_b=\P_{a_1} \cap
\ldots \cap \P_{a_k}$
\item A representative $a\in \mathbb{F}_1$ is said to outer-approximate
(inner-approximate) a representative $b \in\mathbb{F}_2$ if and only if $\P_b
\subseteq \P_a$ ($\P_a \subseteq \P_b$)
\end{itemize}
\end{defn}

\subsubsection{Lower/upper probabilities}

In this paper, lower probabilities (lower
previsions assigned to events) are sufficient to our purpose of representing uncertainty. We
define a {\em lower probability} $\underline{P}$ on $X$ as a
super-additive capacity. Its conjugate
$\overline{P}(A)=1-\underline{P}(A^c)$ is called an upper
probability. The (possibly empty) credal set $\P_{\underline{P}}$ induced by a given
lower probability is its core:
$$\P_{\underline{P}} = \{P \in \setpr{X}| \forall A \subset X, \; P(A) \geq \underline{P}(A)\}. $$
Conversely, from any given non-empty credal set $\P$, one can consider a lower
envelope $P_*$ on events, defined for any event $A \subseteq X$ by
$P_*(A)=\min_{P \in \P} P(A)$. A lower envelope is a super-additive
capacity, and consequently a lower probability. The upper envelope
$P^*(A)=\max_{P \in \P} P(A)$ is the conjugate of $P_*$. In general,
a credal set $\P$ is included in the core of its lower envelope:  $\P \subseteq \P_{P_*}$, since $\P_{P_*}$ can be seen as a
projection of $\P$ on events.

{\em Coherent} lower probabilities $\underline{P}$ are lower
probabilities that coincide with the lower envelopes of their core,
i.e. for all events $A$ of $X$, $\underline{P}(A)=\min_{P \in
\P_{\underline{P}}} P(A).$ Since all representations considered in
this paper correspond to particular instances of
coherent lower probabilities, we will restrict ourselves to such
lower probabilities. In other words, credal sets
$\P_{\underline{P}}$ in this paper are entirely
characterized by their lower probabilities on events and are such
that for every event $A$, there is a probability distribution $P$ in
$\P_{\underline{P}}$ such that $P(A)=\underline{P}(A)$.

A credal set $\P_{\underline{P}}$ can also be described by a set of
constraints on probability assignments to elements of $X$:
\begin{displaymath}
\underline{P}(A) \leq \sum_{x \in A} p(x) \leq \overline{P}(A).
\end{displaymath}
Note that $2^{|X|}-2$ values ($|X|$ being the cardinality of $X$),
are needed in addition to constraints
$\underline{P}(X)=1,\underline{P}(\emptyset)=0$ to completely
specify $\P_{\underline{P}}$.

\subsubsection{Simplified representations}
\label{sec:SimplRepImpProb}  Representing general credal sets induced or not by coherent lower probabilities is usually
costly and dealing with them presents many computational challenges (See, for example,
Walley~\cite{Walley96} or the special issue~\cite{IJARImpProb44}).
In practice, using simpler representations of imprecise
probabilities often alleviates the elicitation and computational burden. 
P-boxes and interval probabilities are two such simplified
representations.

\textbf{P-boxes}

Let us first recall some background on cumulative distributions. Let
$P$ be a probability measure  on the real line $\mathbb{R}$. Its
{\em cumulative distribution}
 is a non-decreasing mapping from $\mathbb{R}$ to $[0,1]$ denoted $F^P$, such that for any $r \in \mathbb{R}$, $F^P(r)=P((-\infty, r])$. Let $F_1$ and $F_2$ be two cumulative distributions. Then, $F_1$ is
said to stochastically dominate $F_2$ if only if $F_1$ is point-wise
lower than $F_2$: $F_1 \leq F_2$.

A p-box~\cite{FersonAll03} is then defined as a pair of (discrete)
cumulative distributions $[\underline{F},\overline{F}]$ such that
$\underline{F}$ stochastically dominates $\overline{F}$. A p-box
induces a credal set $\P_{[\underline{F},\overline{F}]}$ such that:
\begin{equation}\label{eq:usualpbox}
\P_{[\underline{F},\overline{F}]}\!\!=\!\!\{P \in \setpr{\mathbb{R}}
| \forall r \in \mathbb{R}, \; \underline{F}(r) \leq P((-\infty, r])
\leq \overline{F}(r) \}\end{equation} We can already notice that
since sets $(-\infty, r]$ are nested,
$\P_{[\underline{F},\overline{F}]}$ is described by constraints that
are lower and upper bounds on such nested sets (as noticed by Kozine
and Utkin~\cite{KozineUtkin05}, who discuss the problem of building
p-boxes from partial information). This interesting characteristic
will be crucial in the generalized form of p-box we introduce in
section~\ref{sect:GenPboxes}. Conversely we can extract a p-box from
a credal set $\P$ by considering its lower and upper envelopes
restricted to events of the form $(-\infty, r]$, namely, letting  $
\underline{F}(r) = P_*((-\infty, r]),\overline{F}(r) =
P^*((-\infty, r])$. $\P_{[\underline{F},\overline{F}]}$ is then the
tightest outer-approximation of $\P$ induced by a p-box.

Cumulative distributions are often used to elicit probabilistic
knowledge from experts~\cite{Cooke91}. P-boxes can thus directly
benefit from such methods and tools, with the advantages of allowing
some imprecision in the representation (e.g., allowing experts to
give imprecise percentiles). P-boxes are also sufficient to
represent final results produced by  imprecise probability models
when only a threshold violation has to be checked. Working with
p-boxes also allows, via so-called probabilistic
arithmetic~\cite{WilliamsonDowns90}, for very efficient numerical methods
to achieve some particular types of (conservative) inferences.

\textbf{Probability intervals}

 {\em Probability intervals}, extensively studied by De
Campos \emph{et al.}~\cite{CamposAll94},  are defined as
lower and upper bounds of probability distributions. They are
defined by a set of numerical intervals $L=\{ [l(x),u(x)] | x \in X \}$ such
that $l(x) \leq p(x) \leq u(x), \forall x\in X$, where $p(x) =
P(\{x\})$.  A probability interval
induces the following credal set:
\begin{displaymath}
\P_L = \{P \in \setpr{X} | \forall x \in X, \;  l(x) \leq p(x) \leq
u(x) \}
\end{displaymath}

A probability interval $L$ is called reachable if the credal set
$\P_L$ is not empty and if for each element $x \in X$, we can find
at least one probability measure $P \in \P_L$ such that $p(x)=l(x)$
and
 one for which $p(x)=u(x)$. In other words, each bound can be
reached by a probability measure in $\P_L$. Non-emptiness and
reachability respectively correspond to the
conditions~\cite{CamposAll94}:
\begin{align*}
& \sum_{x \in X} l(x) \leq 1 \leq \sum_{x \in X} u(x)  & \textrm{ non-emptiness } \\
&  u(x) + \sum_{y \in X\setminus\{x\}} l(y)  \leq 1 \textrm{ and } l(x) + \sum_{y \in X\setminus\{x\}} u(y)  \geq 1  & \textrm{ reachability } \\
\end{align*}
If a probability interval $L$ is non-reachable,
it can be transformed
 into a probability interval $L'$, by letting  $l'(x)=\inf_{P
\in \P_L} (p(x))$ and $u'(x)= \sup_{P \in \P_L}(p(x))$.
More generally, coherent lower and upper probabilities
$\underline{P}(A),\overline{P}(A)$ induced by $\P_L$ on all events
$A \subset X$ are easily calculated by the following expressions
\begin{equation}
\label{eq:EventProbInt} \underline{P}(A)=\max(\sum_{x \in  A} l(x),
1 - \sum_{x \in A^c} u(x)), \overline{P}(A)=\min(\sum_{x \in  A}
u(x), 1 - \sum_{x \in A^c} l(x)).
\end{equation}
De Campos \emph{et al.}~\cite{CamposAll94} have shown that these
lower and upper probabilities are Choquet capacities of order 2.



Probability intervals, which are modeled by $2|X|$ values, are very
convenient tools to model uncertainty on multinomial data, where
they can express lower and upper confidence bounds. They
can thus be derived from a sample of small size \cite{MassonDenoeux06}.
On the real line, discrete probability intervals correspond to
imprecisely known histograms.  Probability intervals can be
extracted, as useful information, from any credal set $\P$ on a
finite set $X$, by constructing $L_{\P} = \{[\underline{P}(\{x\}),
\overline{P}(\{x\})], x \in X\}$. $L_{\P}$ then represents the
tightest probability interval outer-approximating $\P$.
Numerical and computational advantages that probability
intervals  offer are discussed by De Campos \emph{et
al.}~\cite{CamposAll94}.

\subsection{Random disjunctive sets}

A more specialized setting for representing partial knowledge is
that of random sets. Formally a random set is a family of subsets of $X$ each bearing a probability weight.  Typically, each set represents an incomplete
observation, and the probability bearing on this set should
potentially be shared among its elements, but is not by lack of
sufficient information.

\subsubsection{Belief and Plausibility functions}
Formally, a random set is defined as a mapping $\Gamma:\Omega\to\wp(X)$ from a
probability space $(\Omega, \mathcal{A}, P)$ to the power set
$\wp(X)$ of another space $X$ (here finite). It is also called a
multi-valued mapping $\Gamma$. Insofar as sets $\Gamma(\omega)$
represent incomplete knowledge about a single-valued random variable, each such set
contains mutually exclusive elements and is called a {\em disjunctive}
set\footnote{as opposed to sets as collections of objects, i.e.
sets whose elements are jointly present, such as a region in a
digital image.}. Then this mapping induces the following coherent
lower and upper probabilities on $X$ for all events
$A$~\cite{Dempster67} (representing all probability functions on X
that could be found if the available information were complete):
\begin{displaymath}
\underline{P}(A)=P(\{\omega \in \Omega | \Gamma(\omega) \subseteq A
\})  ; \quad \overline{P}(A)=P(\{\omega \in \Omega | \Gamma(\omega)
\cap A \neq \emptyset \}),
\end{displaymath}
where $\{\omega \in \Omega | \Gamma(\omega)
\cap A \neq \emptyset \} \in \mathcal{A}$ is assumed. When $X$ is finite, a random set can be represented as a mass
assignment $m$ over the power set $\wp(X)$ of $X$, letting $m(E) =
P(\{\omega,\Gamma(\omega) = E\}), \forall E \in X$. Then, $\sum_{E
\subseteq X} m(E) = 1$ and $m(\emptyset) = 0$. A set $E$ that
receives strict positive mass is called a focal set, and the mass
$m(E)$ can be interpreted as the probability that the most precise
description of what is known about a particular situation is of the
form "$x \in E$". From this mass assignment, Shafer~\cite{Shafer76}
define two set functions, called {\em belief and plausibility
functions}, respectively:
\begin{displaymath}
Bel(A) = \sum_{E,E \subseteq A} m(E); \qquad Pl(A) = 1 - Bel(A^c) =
\sum_{E, E \cap A \neq \emptyset} m(E).
\end{displaymath}
The mass assignment being positive, a belief function is an
$\infty$-monotone capacity. The mass assignment $m$
 is indeed the M\"obius transform of the capacity $Bel$. Conversely, any
$\infty$-monotone capacity is induced by one and only one random
set. We can thus speak of {\em the} random set underlying $Bel$. In the
sequel, we will use this notation for lower probabilities stemming
from random sets (Dempster and Shafer definitions being equivalent
on finite spaces). Smets~\cite{Smets05} has studied the case of
continuous random intervals defined on the real line $\mathbb{R}$,
where the mass function is replaced by a mass density bearing on
pairs of interval endpoints.

Belief functions can be considered as special cases of coherent
lower probabilities, since they are $\infty$-monotone capacities. A
random set thus induces the credal set $\P_{Bel}= \{ P \in \setpr{X}
| \forall A \subseteq X,  Bel(A) \leq P(A) \}$.

Note that Shafer~\cite{Shafer76} does not refer to an underlying
probability space, nor does he uses the fact that a belief function
is a lower probability: in his view, extensively taken over by
Smets~\cite{Smets97}, $Bel(A)$ is supposed to quantify an agent's
belief per se with no reference to a probability. However, the
primary mathematical tool common to Dempster's upper and lower
probabilities and to the Shafer-Smets view is the notion of
(generally finite) random disjunctive set.

\subsubsection{Practical aspects}

In general, $2^{|X|}-2$ values are still needed to completely
specify a random set, thus not clearly reducing the complexity of
the model representation with respect to capacities. However, simple
belief functions defined by only a few positive focal elements do
not exhibit such complexity.
For instance, a simple support belief
function is a natural model of an unreliable testimony, namely an
expert stating that the value of a parameter $x$ belong to set  $A
\subseteq X$. Let $\alpha$ be the reliability of the expert
testimony, i.e. the probability that the information is irrelevant.
The corresponding mass assignment is $m(A) = \alpha, m(X) = 1-
\alpha$.
 Imprecise results from some statistical experiments
are easily expressed by means of random sets, $m(A)$ being the
probability of an observation of the form $x \in A$.

As practical models of uncertainty, random sets present many
advantages. First, as they can be seen as probability distributions
over subsets of $X$, they can be easily simulated by classical
methods such as Monte-Carlo sampling, which is not the case for
other Choquet capacities. On the real line, a random set is often
restricted to a finite collection of closed intervals with
associated weights, and one can then easily extend results from
interval analysis~\cite{Moore79} to random
intervals~\cite{DuboisPrade91,OberkampfAll04Bis}.


\subsection{Possibility theory}

The primary mathematical tool of possibility theory is the
possibility distribution, which is a set-valued piece of information
where some elements are more plausible than others.  To a
possibility distribution  are associated specific measures of
certainty and plausibility.

\subsubsection{Possibility and necessity measures} A possibility
distribution is a mapping $\pi:X \to [0,1]$ from $X$ to the unit
interval such that $\pi(x) = 1$ for at least one element $x$ in $X$.
Formally, a possibility distribution is equivalent to the membership
function of a normalized fuzzy set~\cite{Zadeh78}\footnote{The membership
function of a fuzzy set $\nu$ is a mapping $\nu:X\to [0,1]$}. Twenty
years earlier, Shackle~\cite{shackle1961} had introduced an
equivalent notion called distribution of potential surprise
(corresponding to $1 - \pi(x)$) for the representation of
non-probabilistic uncertainty.

Several set-functions can be defined from a possibility distribution
$\pi$~\cite{DuboisAll00}:
\begin{eqnarray} \label{eq:possnecdegree}
\mbox{Possibility measures: } & \Pi(A)  & = \sup_{x\in A} \pi(x). \\
\mbox{Necessity measures: }   & N(A)  & =1 - \Pi(A^c).\\
\mbox{Sufficiency measures: } & \Delta(A) & = \inf_{x\in A} \pi(x).
\end{eqnarray}
The possibility degree of an event $A$ evaluates the extent to which
this event is plausible, i.e., consistent with the available
information. Necessity degrees express the certainty of events, by
duality. In this context, distribution $\pi$ is potential (in the
spirit of Shackle's), i.e. $\pi(x) = 1$ does not guarantee the
existence of $x$. Their characteristic property are: $N(A\cap B) =
\min(N(A), N(B))$  and  $\Pi(A\cup B) = \max(\Pi(A), \Pi(B))$ for
any pair of events $A,B$ of $X$. On the contrary $\Delta(A)$
measures the extent to which all states of the world where $A$
occurs are plausible. Sufficiency\footnote{also called guaranteed possibility distributions~\cite{DuboisAll00}.}
distributions, generally denoted by $\delta$, express {\em actual} possibility. They are
understood as degree of empirical support and obey an opposite
convention: $\delta(x)$ = 1 guarantees (= is sufficient for) the
existence of $x$.

\subsubsection{Relationships with previous theories}



A necessity measure  (resp a possibility measure) is formally a particular case of
belief function (resp. a plausibility function) induced by a random
set with nested focal sets (already in \cite{Shafer76}). Given a possibility distribution $\pi$ and a degree $\alpha \in
[0,1]$, strong and regular $\alpha$-cuts are subsets respectively
defined as \mbox{$A_{\overline{\alpha}} = \{x \in X | \pi(x) >
\alpha \}$} and $A_{\alpha} = \{x \in X | \pi(x) \geq \alpha \}$.
These $\alpha$-cuts are nested, since if $\alpha > \beta$, then
$A_{\alpha} \subseteq A_{\beta}$. On finite spaces, the set  $\{
\pi(x), x\in  X\}$ is of the form  \mbox{$\alpha_0 = 0 < \alpha_1 <
\ldots < \alpha_M=1$}. There are only  $M$ distinct
$\alpha$-cuts. A
possibility distribution $\pi$ then induces a random set having, for
$i=1,\ldots,M$, the following focal sets $E_i$ with masses $m(E_i)$:
\begin{equation}
\label{eq:possRStransform} \left\{ \begin{array}{c} E_i=\{ x \in X
| \pi(x) \geq \alpha_i \} = A_{\alpha_i} \\
m(E_i)=\alpha_i - \alpha_{i-1}
\end{array} \right.
\end{equation}
In this nested situation, the same amount of information is
contained in the mass function $m$ and the possibility distribution
$\pi(x) = Pl(\{x\})$, called the {\em contour} function of $m$. For instance a simple support belief function
such that $ m(A) = \alpha, m(X) = 1 - \alpha$ forms a nested
structure, and yields the possibility distribution $\pi(x) = 1$ if
$x \in A$ and $1-\alpha$ otherwise. In the general case, $m$ cannot
be reconstructed only from its contour function. Outer and inner approximations of
general random sets in terms of possibility distributions have been
studied by Dubois and Prade in~\cite{DuboisPrade90}.

Since the necessity measure is formally a particular case of belief function, it is
also an $\infty$-monotone capacity, hence  a particular coherent
lower probability. If the necessity measure is viewed as a coherent
lower probability, its possibility distribution induces the credal
set $\P_{\pi}=\{P \in \setpr{X}| \forall A \subseteq X, \; P(A) \geq
N(A)\}$. We recall here a result, proved by Dubois
\emph{et al.}~\cite{DuboisPrade92bis,DuboisAll04} and by Couso \emph{et
al.}~\cite{CousoAll01} in a more general setting, which links
probabilities $P$ that are in $\P_{\pi}$ with constraints on
$\alpha$-cuts, that will be useful in the sequel:
\begin{prop}
\label{prop:LinkAlphCutImprProb} Given a possibility distribution
$\pi$ and the induced convex set $\P_{\pi}$, we have for all
$\alpha$ in $(0, 1]$, $P \in \P_{\pi}$ if and only if $$1- \alpha
\leq P( \{ x \in X | \pi(x)
> \alpha \})$$
\end{prop}
This result means that the probabilities $P$ in the credal set
$\P_{\pi}$ can also be described in terms of constraints on strong
$\alpha$-cuts of $\pi$ (i.e. $1- \alpha \leq P(
A_{\overline{\alpha}})$).

\subsubsection{Practical aspects}

At most $|X|-1$ values are needed to fully assess a possibility
distribution, which makes it the simplest uncertainty representation
explicitly coping with imprecise or incomplete knowledge. This
simplicity makes this representation very easy to handle. This also
implies less expressive power, in the sense that, for any event $A$
, either $\Pi(A) = 1$  or $N(A) = 0$  (i.e. intervals
$[N(A),\Pi(A)]$ are of the form $[0,\alpha]$ or $[\beta,1]$). This
means that, in several situations, possibility distributions will be
insufficient to reflect the available information.

Nevertheless, the expressive power of possibility distributions fits
various practical situations.  Moreover, a recent
psychological study~\cite{RaufasteAll03} shows that sometimes people handle uncertainty according to possibility theory rules. Possibility distributions on the real
line can be interpreted as a
set of nested intervals with different confidence degrees~\cite{DuboisAll04} (the larger the
set, the highest the confidence degree), which is a good model of, for example, an expert
opinion concerning the value of a badly known parameter. Similarly, it is
natural to view nested confidence intervals coming from statistics
as a possibility distribution. Another practical case where
uncertainty can be modeled by possibility distributions is the case
of vague linguistic assessments concerning
probabilities~\cite{DeCoo05}.

\subsection{P-boxes and probability intervals in the uncertainty landscape}

P-boxes, reachable probability intervals, random sets and
possibility distributions can all be modeled by credal sets and
define coherent lower probabilities.  Kriegler and
Held~\cite{KrieglerHeld05} show that random sets generalize p-boxes
(in the sense of Definition~\ref{def:compa}), but that the converse
do not hold (i.e. credal sets induced by different random sets can
have the same upper and lower bounds on events of the type $(\infty,r]$, and hence induce the same p-boxes).

There is no specific relationship between the frameworks of possibility
distributions, p-boxes and probability intervals, in the sense that
none generalize the other. Some results comparing possibility
distributions and p-boxes are given by Baudrit and
Dubois~\cite{BaudritDubois06}. Similarly, there is no generalization
relationship between probability intervals and random sets. Indeed
upper and lower probabilities induced by reachable probability
intervals are order 2 capacities only, while belief functions are
$\infty$-monotone.  In general, one can only approximate one
representation by the other.





Transforming a belief function $Bel$ into the tightest probability interval $L$ outer-approximating it (i.e. $\P_{Bel} \subset
\P_L$, following Definition~\ref{def:compa}) is simple, and consists
of taking for all $x \in X$:
\begin{displaymath}
l(x) = Bel(\{x\}) \textrm{ and } u(x) = Pl(\{x\})
\end{displaymath}
and since belief and plausibility functions are the lower envelope
of the induced credal set $\P_{Bel}$, we are sure that the so-built
probability interval $L$ is reachable.


The converse problem, i.e. to  transform a set $L$ of probability
intervals into an inner-approximating random set was studied by
Lemmer and Kyburg~\cite{LemmerKyburg91}. On the contrary,
Denoeux~\cite{Denoeux06} extensively studies the problem of
transforming a probability interval $L$ into a random set
outer-approximating $L$ (i.e., $\P_L \subset \P_{Bel}$). The
transformation of a given probability interval $L$ into an
outer-approximating possibility distribution is studied by Masson
and Denoeux~\cite{MassonDenoeux06}, who propose efficient methods to
achieve such a transformation.

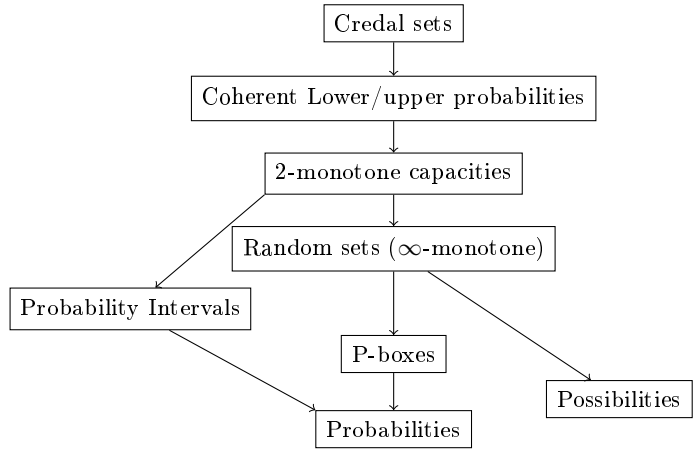
\begin{figure}
\begin{center}
\begin{tikzpicture}
 \node (ImpProb) at (0,-0.4) [draw] {\scriptsize{Credal sets}};
 \node (LowUpp) at (0,-1.4) [draw] {\scriptsize{Coherent Lower/upper probabilities}};
 \node (2mon) at (0,-2.4) [draw] {\scriptsize{2-monotone capacities}} ;
 \node (randset) at (0,-3.4) [draw] {\scriptsize{Random sets ($\infty$-monotone)}} ;
 \node (Pbox) at (0,-4.8) [draw] {\scriptsize{P-boxes}} ;
 \node (Prob) at (0,-5.8) [draw] {\scriptsize{Probabilities}} ;
 \node (ProbInt) at (-3.5,-4.2) [draw] {\scriptsize{Probability
 Intervals}} ;
 \node (Poss) at (3,-5.4) [draw] {\scriptsize{Possibilities}} ;

 \draw[->] (ImpProb) -- (LowUpp)  ;
 \draw[->]  (LowUpp) -- (2mon) ;
 \draw[->]  (2mon) -- (randset)  ;
 \draw[->]  (randset) -- (Pbox) ;
 \draw[->]  (Pbox) -- (Prob) ;
 \draw[->] (randset) -- (Poss) ;
 \draw[->] (2mon.south west) -- (ProbInt)  ;
 \draw[->]  (ProbInt) -- (Prob.north west) ;

\end{tikzpicture}  \caption{Representation relationships: summary  $A \longrightarrow B$: $A$ generalizes $B$}\label{fig:summarystate}
\end{center}
\end{figure}

The main relations existing between imprecise probabilities,
lower/upper probabilities, random sets, probability intervals,
p-boxes and possibility distributions, are pictured on
Figure~\ref{fig:summarystate}. From top to bottom, it goes from the
more general, expressive and complex theories to the less general,
less expressive but simpler representations. Arrows are directed
from a given representation to the representations it generalizes.

\section{Generalized p-boxes}
 \label{sect:GenPboxes}

As recalled in Section~\ref{sec:ImpProb},  p-boxes are useful
representations of uncertainty in many practical applications\cite{CooperAll96,FersonAll03,KrieglerHeld05}. So
far, they only make sense on the (discretized) real line equipped
with the natural ordering of numbers. P-boxes are instrumental to extract interpretable information from imprecise probability representations. They provide faithful estimations of the probability that a variable $\tilde{x}$ violates a threshold $\theta$, i.e.,  upper and lower estimates of the probability of events of the form $\tilde{x} \geq \theta$. However, they are much less adequate to compute the probability that some output remains close to a reference value $\rho$ \cite{BaudritDubois06}, which corresponds to computing upper and lower estimates of the probability of events of the form $\vert \tilde{x} - \rho \vert \geq \theta$.  
The rest of the paper is devoted to the study of  a generalization 
p-boxes, to arbitrary (finite) spaces, where the underlying ordering relation is arbitrary, and that can address this type of query. Generalized p-boxes will also be
instrumental to characterize a recent representation proposed by
Neumaier~\cite{Neumaier04}, studied in the second part of this paper.

Generalized p-boxes are defined in Section~\ref{sec:GenPboxDef}. We
then proceed to show the link between generalized p-boxes,
possibility distributions and random sets. We first show  that the former generalize
possibility distributions and are representable (in the sense of Definition~\ref{def:compa}) by pairs thereof. 
Connections between generalized p-boxes and probability intervals are also explored.

%
%
%
%

\subsection{Definition of generalized p-boxes}
\label{sec:GenPboxDef}

The starting point of our generalization is to notice that any two cumulative distribution functions modelling a p-box are comonotonic. Two mappings $f$ and $f'$ from a space $X$ to the real line are said
to be comonotonic if and only if, for any pair of elements $x,y \in
X$, we have $f(x) < f(y) \Rightarrow f'(x) \leq f'(y)$. In other
words, given an indexing of $X=\{x_1,\ldots,x_n\}$, there is a
permutation $\sigma$ of $\{1, 2, \dots, n\}$ such that
$f(x_{\sigma(1)})\geq f(x_{\sigma(2)})\geq \dots \geq f(x_{\sigma(n)})$
and $f'(x_{\sigma(1)})\geq f'(x_{\sigma(2)})\geq \dots \geq
f'(x_{\sigma(n)})$. We define a generalized p-box as follows:
\begin{defn}
\label{def:genpboxes} A generalized p-box
$[\underline{F},\overline{F}]$ defined on $X$ is a pair of
comonotonic mappings $\underline{F},\overline{F}$, $\underline{F}: X
\to [0,1]$ and $\overline{F}: X \to [0,1]$ from $X$ to $[0,1]$ such
that $\underline{F}$ is pointwise less than $\overline{F}$ (i.e.
$\underline{F} \leq \overline{F} $) and there is at least one
element $x$ in $X$ for which $\overline{F}(x)=\underline{F}(x)=1$.
\end{defn}
Since each distribution $\underline{F},\overline{F}$ is fully
specified by $|X| -1$ values, it follows that $2|X|-2$ values
completely determine a generalized p-box. Note that, given a
generalized p-box $[\underline{F},\overline{F}]$, we can always
define a {\em complete} pre-ordering
$\leq_{[\underline{F},\overline{F}]}$ on  $X$ such that $x
\leq_{[\underline{F},\overline{F}]} y$ if $\underline{F}(x) \leq
\underline{F}(y)$ and \mbox{$\overline{F}(x) \leq \overline{F}(y)$},
due to the comonotonicity condition. If $X$ is a subset of the real
line and if $\leq_{[\underline{F},\overline{F}]}$ is the natural
ordering of numbers, then we retrieve classical p-boxes.

To simplify notations in the sequel, we will consider that, given a
generalized p-box $[\underline{F},\overline{F}]$, elements $x$ of
$X$ are indexed such that $i < j$ implies that $x_i
\leq_{[\underline{F},\overline{F}]} x_j$. We will denote
$(x]_{[\underline{F},\overline{F}]}$ the set of the form $\{x_i |
x_i \leq_{[\underline{F},\overline{F}]} x \}$. The credal set
induced by a generalized p-box $[\underline{F},\overline{F}]$ can
then be defined as
$$\P_{[\underline{F},\overline{F}]}=\{P \in \setpr{X} |i=1,\ldots,n, \; \underline{F}(x_i) \leq P((x_i]_{[\underline{F},\overline{F}]}) \leq \overline{F}(x_i) \}.$$
It induces coherent upper and lower probabilities such that
$\underline{F}(x_i)=
\underline{P}((x_i]_{[\underline{F},\overline{F}]})$ and
$\overline{F}(x_i)=\overline{P}((x_i]_{[\underline{F},\overline{F}]})$.
When $X=\mathbb{R}$ and $\leq_{[\underline{F},\overline{F}]}$ is the
natural ordering on numbers, then $\forall r \in \mathbb{R}$,
$(r]_{[\underline{F},\overline{F}]}=(-\infty,r]$, and the above
equation coincides with Equation (\ref{eq:usualpbox}).

In the following, sets $(x_i]_{[\underline{F},\overline{F}]}$ are
denoted $A_i$, for all $i=1,\ldots,n$. These sets are nested, since
 $\emptyset \subset A_1 \subseteq \ldots \subseteq A_n =
X$\footnote{Since $\leq_{[\underline{F},\overline{F}]}$ is a complete pre-order on $X$, we can have
$x_i =_{[\underline{F},\overline{F}]} x_{i+1}$ and $A_i=A_{i+1}$,
which explains the non-strict inclusions. They would be strict if
$<_{[\underline{F},\overline{F}]}$ were a linear order.}. For all
$i=1,\ldots,n$, let $\underline{F}(x_i)=\alpha_i$ and
$\overline{F}(x_i)=\beta_i$. With these conventions, the credal set
$\P_{[\underline{F},\overline{F}]}$ can now be described by the
following $n$ constraints bearing on probabilities of nested sets
$A_i$:
\begin{equation}
\label{eq:PBoxCons} i=1,\ldots,n \qquad \alpha_i \leq P(A_i) \leq
\beta_i
\end{equation}
with $0 = \alpha_0 \leq \alpha_1 \leq \ldots \leq \alpha_n = 1$, $0
= \beta_0 < \beta_1 \leq \beta_2 \leq \ldots \leq \beta_n = 1$ and
$\alpha_i \leq \beta_i$.


As a consequence, a generalized p-box can be generated in two
different ways:
\begin{itemize}
\item By starting from two comonotone functions $\underline{F}\leq \overline{F}$
 defined on $X$, the pre-order being induced by the values of these functions,
\item or by assigning upper and lower bounds on probabilities of a prescribed collection of nested
sets $A_i$.
\end{itemize}
Note that the second approach is likely to be more useful in
practical assessments and elicitation of generalized p-boxes.

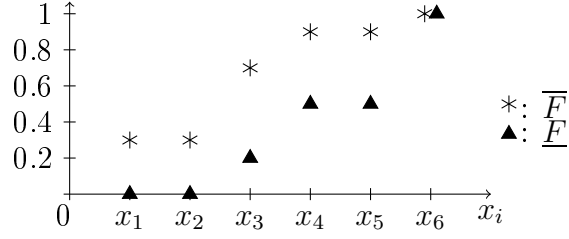
\begin{figure}
\begin{center}
\begin{tikzpicture}[scale=0.8]
\draw[->] (-0.1,0) node[below] {0} -- (7,0) node[below] {$x_i$} ;
\draw[->] (0,-0.1) -- (0,3.2) ; \draw (0.1,3) -- (-0.1,3) node[left]
{1} ; \draw (1,0.1) --  (1,-0.1) node[below] {$x_1$}; \draw (2,0.1)
-- (2,-0.1) node[below] {$x_2$}; \draw (3,0.1) -- (3,-0.1)
node[below] {$x_3$}; \draw (4,0.1) --  (4,-0.1) node[below] {$x_4$};
\draw (5,0.1) --  (5,-0.1) node[below] {$x_5$}; \draw (6,0.1) --
(6,-0.1) node[below] {$x_6$}; \draw (0.1,0.6) -- (-0.1,0.6)
node[left] {0.2} ; \draw (0.1,1.2) -- (-0.1,1.2) node[left] {0.4} ;
\draw (0.1,1.8) -- (-0.1,1.8) node[left] {0.6} ; \draw (0.1,2.4) --
(-0.1,2.4) node[left] {0.8} ;

\draw plot[only marks,mark=triangle*,mark options={scale=2.0}]
coordinates {(1,0) (2,0) (3,0.6) (4,1.5) (5,1.5) (6.1,3)} ;

\draw plot[only marks,mark=asterisk,mark options={scale=2.0}]
coordinates {(1,0.9) (2,0.9) (3,2.1) (4,2.7) (5,2.7) (5.9,3)};

\draw plot[mark=asterisk,mark options={scale=2.0}] coordinates
{(7.3,1.5)} ; \node[right] at (7.3,1.5) {: $\overline{F}$} ;

\draw plot[mark=triangle*,mark options={scale=2.0}] coordinates
{(7.3,1)} ; \node[right] at (7.3,1) {: $\underline{F}$} ;

\end{tikzpicture}
\caption{Generalized p-box $[\underline{F},\overline{F}]$ of
Example~\ref{exmp:genpbox1}} \label{fig:PboxEx}
\end{center}
\end{figure}

\begin{exmp}
\label{exmp:genpbox1} All along this section, we will use this
example to illustrate results on generalized p-boxes. Let $X = \{ x_1,\ldots,x_6 \}$. These
elements could be, for instance, the facets of a biased die.
For various reasons, we only have incomplete information about
the probability of some subsets $A_1 = \{ x_1, x_2 \}$, $A_2 = \{
x_1,x_2,x_3 \}$, \mbox{$A_3 = \{ x_1,\ldots, x_5 \}$}, or $X(=A_4)$.
An expert supplies the
following confidence bounds on the  frequencies of these sets:
\begin{displaymath}
P(A_1) \in  [0,0.3] \qquad P(A_2) \in [0.2,0.7] \qquad P(A_3) \in
[0.5,0.9]
\end{displaymath}
The uncertainty can be modeled by the generalized p-box pictured on
Figure~\ref{fig:PboxEx}:
\begin{displaymath}
\begin{array}{c@{\hspace{10mm}}c@{\hspace{4mm}}c@{\hspace{4mm}}c@{\hspace{4mm}}c@{\hspace{4mm}}c@{\hspace{4mm}}c}
& x_1  & x_2 & x_3 & x_4 & x_5 & x_6 \\
\hline
\overline{F} & 0.3  & 0.3 & 0.7 & 0.9 & 0.9 & 1 \\
\underline{F} & 0  & 0 & 0.2 & 0.5 & 0.5 & 1. \\
\hline
\end{array}
\end{displaymath}
\end{exmp}

\subsection{Connecting generalized p-boxes with possibility distributions}
\label{sec:GenpboxPoss}


It is natural to search for a connection between generallized p-boxes
and possibility theory, since possibility distributions can be
interpreted as a collection of nested sets with associated lower
bounds, while generalized p-boxes correspond to lower and upper
bounds also given on a collection of nested sets. Given a
generalized p-box $[\underline{F},\overline{F}]$, the following
proposition holds:
\begin{prop}
\label{prop:PBoxPoss} Any generalized p-box
$[\underline{F},\overline{F}]$ on $X$ is representable by a pair of
possibility distributions $\pi_{\overline{F}},\pi_{\underline{F}}$,
such that $\P_{[\underline{F},\overline{F}]} =
\P_{\pi_{\overline{F}}}\cap \P_{\pi_{\underline{F}}}$, where:
$$\pi_{\overline{F}}(x_i)=\beta_i \mbox{~and~}
\pi_{\underline{F}}(x_i)= 1 - \max\{\alpha_j | 
\alpha_j < \alpha_{i}; j=0,\ldots,i \}$$ for $i=1,\ldots,n$, with $\alpha_0=0$.
\end{prop}
\begin{proof}[\textbf{Proof of Proposition~\ref{prop:PBoxPoss}}]
Consider the set of constraints given by Equation
(\ref{eq:PBoxCons}) and defining the convex set
$\P_{[\underline{F},\overline{F}]}$. These constraints can be
separated into two distinct sets: $(P(A_i) \leq \beta_i)_{i= 1, n}$
and $(P(A_i^c) \leq 1-\alpha_i)_{i= 1, n}$. Now, rewrite constraints
of Proposition~\ref{prop:LinkAlphCutImprProb},  in the form $\forall
\alpha \in (0, 1]$: $P \in  \P_{\pi} \textrm{ if and only if }
 P( \{ x \in X | \pi(x) \leq \alpha \}) \leq \alpha $.

The first set of constraints $(P(A_i) \leq \beta_i)_{i= 1,n}$
defines a credal set $\P_{\pi_{\overline{F}}}$ that is induced by
the possibility distribution $\pi_{\overline{F}}$, while the second
set of constraints $(P(A_i^c) \leq 1-\alpha_i)_{i= 1, n}$ defines a
credal set $\P_{\pi_{\underline{F}}}$ that is induced by the
possibility distribution $\pi_{\underline{F}}$, since $A_i^c =
\{x_k, \dots x_n\}$, where $k =  \max\{j | \alpha_j < \alpha_{i}
\}$.
The credal set of the generalized p-box
$[\underline{F},\overline{F}]$, resulting from the two sets of
constraints, namely  $i=1,\ldots,n, \quad \beta_i \leq P(A_i) \leq
\alpha_i$, is thus $\P_{\pi_{\overline{F}}} \cap
\P_{\pi_{\underline{F}}}$.
\end{proof}

\begin{exmp}
\label{exmp:genpbox2} The possibility distributions
$\pi_{\overline{F}},\pi_{\underline{F}}$ for the generalized p-box
defined in Example~\ref{exmp:genpbox1} are:
\begin{displaymath}
\begin{array}{c@{\hspace{10mm}}c@{\hspace{4mm}}c@{\hspace{4mm}}c@{\hspace{4mm}}c@{\hspace{4mm}}c@{\hspace{4mm}}c}
& x_1  & x_2 & x_3 & x_4 & x_5 & x_6 \\
\hline
\pi_{\overline{F}} & 0.3  & 0.3 & 0.7 & 0.9 & 0.9 & 1 \\
\pi_{\underline{F}} & 1 & 1 & 1 & 0.8 & 0.8 & 0.5 \\
\hline
\end{array}
\end{displaymath}
\end{exmp}

Note that, when $\underline{F}$ is injective,
$<_{[\underline{F},\overline{F}]}$ is a linear order, and we have
$\pi_{\underline{F}}(x_i)= 1 - \alpha_{i-1}$. So, generalized
p-boxes allow to model uncertainty in terms of pairs of comonotone
possibility distributions. In this case, contrary to the case of
only one possibility distribution, the two bounds enclosing
the probability of a particular event $A$ can be tighter, i.e. no longer
restricted to the form $[0,\alpha]$ or $[\beta,1]$, but
contained in the intersection of intervals of this form.

An interesting case is the one where, for all $i=1,\ldots,n-1$,
$\underline{F}(x_i)=0$ and $\underline{F}(x_n)=1$. Then,
$\pi_{\underline{F}} = 1$ and $\P_{\pi_{\overline{F}}} \cap
\P_{\pi_{\underline{F}}} =\P_{\pi_{\overline{F}}}$ and we retrieve
the single distribution $\pi_{\overline{F}}$. We recover
$\pi_{\underline{F}}$ if we take for all $i=1,\ldots,n$,
$\overline{F}(x_i)=1$. This means that generalized p-boxes
also generalize possibility distributions, and are representable by them in the sense of Definition~\ref{def:compa}.

\subsection{Connecting Generalized p-boxes and random sets}
\label{sec:GenPboxRS}

We already mentioned that p-boxes are special cases of random sets, and the
following proposition shows that it is also true for generalized
p-boxes.
\begin{prop}
\label{prop:PBoxBel} Generalized p-boxes are special cases of random sets, in
the sense that for any generalized p-box
$[\underline{F},\overline{F}]$ on $X$, there exist a belief function $Bel$ such that $\P_{[\underline{F},\overline{F}]}= \P_{Bel}$.
\end{prop}
In order to prove Proposition~\ref{prop:PBoxBel}, we show that the
lower probabilities induced by a generalized p-box and by
the belief function given by Algorithm~\ref{alg:PboxBel} coincide on
every event. To do that, we use the partition of $X$ induced by
nested sets $A_i$, and compute lower probabilities of elements of this partition. We then check that the lower
probabilities on all events induced by the generalized p-box
coincide with the belief function induced by
Algorithm~\ref{alg:PboxBel}. The detailed proof can be found in the
appendix.

Algorithm~\ref{alg:PboxBel} below provides an easy way to build the
random set encoding a given generalized p-box. It is similar to
existing algorithms~\cite{KrieglerHeld05,ReganAll04}, and extends
them to more general spaces. The main idea of the algorithm is to
use the fact that a generalized p-box can be seen as a random set
whose focal sets are obtained by thresholding the cumulative
distributions (as in Figure~\ref{fig:PboxEx}). Since  the sets $A_i$
are nested, they induce a partition of $X$ whose elements are of the
form $G_i = A_{i} \setminus A_{i-1}$. The focal sets of the random
set equivalent to the generalized p-box are made of unions of
consecutive elements of this partition. Basically, the procedure
comes down to progressing a threshold $\theta\in [0, 1]$. When
$\alpha_{i+1} > \theta \geq \alpha_i$ and $\beta_{j+1} > \theta \geq
\beta_j$, then, the corresponding focal set is $A_{i+1} \setminus
A_j$, with mass
\begin{equation}
\label{eq:PboxRSpart} m(A_{i+1} \setminus A_j) =
\min(\alpha_{i+1},\beta_{j+1})-
\max(\alpha_i,\beta_j).\end{equation}
\begin{algorithm2e} 
\dontprintsemicolon \KwIn{Generalized p-box $[\underline{F},\overline{F}]$ and corresponding nested sets
$\emptyset=A_0,A_1,\ldots,A_n=X$, lower bounds $\alpha_i$ and upper
bounds $\beta_i$ } \KwOut{Equivalent random set}
\For{$i=1,\ldots,n$}{Build partition
 $ G_i = A_i \setminus A_{i-1} $ } \; Build set $\{\gamma_l | l=1,\ldots,2n-1\} = \{\alpha_i | i=1,\ldots,n \} \cup \{\beta_i | i=1,\ldots,n
 \}$ \;
 With $\gamma_l$ indexed such that $\gamma_1 \leq \ldots \leq \gamma_l \leq \ldots
\leq \gamma_{2n-1}=\beta_n=\alpha_n=1$ \; Set $\alpha_0 = \beta_0 =
\gamma_0 = 0$ and focal set $E_0 = \emptyset$ \;
\For{$k=1,\ldots,2n-1$}{\If{$\gamma_{k-1} = \alpha_i$}{$E_k =
E_{k-1} \ \cup \ G_{i+1}$} \; \If{$\gamma_{k-1} = \beta_i$}{$E_k =
E_{k-1} \ \setminus \ G_{i}$} \; Set $m(E_k)=\gamma_k- \gamma_{k-1}$
} \caption{R-P-box $\rightarrow$ random set transformation}
\label{alg:PboxBel}
\end{algorithm2e}

We can also give another characterization of the random set
(\ref{eq:PboxRSpart}): let us note $0 = \gamma_0 < \gamma_1 < \ldots
< \gamma_M =1$ the distinct values taken by
$\underline{F},\overline{F}$ over elements $x_i$ of $X$ (note that
$M$ is finite and $M < 2n$). Then, for $j=1,\ldots,M$, the random
set defined as:
\begin{equation}
\label{eq:PboxRStransform} \left\{ \begin{array}{c} \scriptstyle
E_j=\{ x_i \in X
| (\pi_{\overline{F}}(x_i) \geq \gamma_j) \wedge (1-\pi_{\underline{F}}(x_i) < \gamma_j ) \} \\
m(E_j)=\gamma_j - \gamma_{j-1}
\end{array} \right.
\end{equation}
is the same as the one built by using Algorithm~\ref{alg:PboxBel},
but this formulation lays bare the link between Equation
(\ref{eq:possRStransform}) and the possibility distributions
$\pi_{\overline{F}},\pi_{\underline{F}}$.

\begin{exmp}
\label{exmp:genpbox3} This example illustrates the application of
Algorithm~\ref{alg:PboxBel}, by applying it to the generalized p-box
given in Example~\ref{exmp:genpbox1}. We have:
\begin{align*}
& G_1 = \{ x_1, x_2 \} \qquad G_2 = \{ x_3 \} \qquad G_3 = \{
x_4,x_5 \} \qquad G_4 = \{ x_6 \} &
\end{align*}
and
\begin{align*}
& 0  \leq 0  \leq 0.2  \leq 0.3  \leq 0.5  \leq 0.7  \leq 0.9  \leq 1  & \\
& \alpha_0 \leq  \alpha_1 \leq \alpha_2 \leq \beta_1 \leq \alpha_3 \leq \beta_2 \leq \beta_3 \leq \alpha_4& \\
& \gamma_0 \leq \gamma_1 \leq \gamma_2 \leq \gamma_3 \leq \gamma_4 \leq \gamma_5 \leq \gamma_6 \leq \gamma_7 & \\
\end{align*}
which finally yields the following random set
\begin{align*}
 & & m(E_1) = m(G_1) = 0 \qquad m(E_2) = m(G_1 \cup G_2) = 0.2 \\
 & & m(E_3) = m(G_1 \cup G_2 \cup G_3) = 0.1 \qquad  m(E_4) = m(G_2 \cup G_3)
= 0.2  \\ & & m(E_5) = m(G_2 \cup G_3 \cup G_4) = 0.2 \qquad m(E_6)
= m(G_3 \cup G_4) = 0.2 \\
& &  m(E_7) = m(G_4) = 0.1 
\end{align*}
This random set can then be used as an alternative representation of
the provided information.
\end{exmp}
Propositions~\ref{prop:PBoxBel} and~\ref{prop:PBoxPoss} together
indicate that generalized p-boxes are more expressive than single
possibility distributions and less expressive than random sets, but, as
recalled before, less expressive (and, in this sense, simpler)
models are often easier to handle in practice. The following explicit expression for lower probabilities induced  by generalized p-boxes $[\underline{F},\overline{F}]$ on $X$
shows that we can expect it to be the case (see appendix):
 \begin{equation}
\label{eq:LowerprobGenPboxes}
\underline{P}(\bigcup_{k=i}^j
G_k)=\max(0,\alpha_j-\beta_{i-1}).
\end{equation}
Let us call a subset $E$ of $X$ \emph{$[\underline{F},\overline{F}]$-connected} if it can expressed as
an union of consecutive elements $G_k$, i.e. $E=\bigcup_{k=i}^j
G_k$, with $0<i<j\leq n$. For any event $A$, let \mbox{$A_*=\bigcup_{ E \subseteq A} E$},
with $E$ all maximal \emph{$[\underline{F},\overline{F}]$-connected} subsets included in $A$. We know  (see appendix)
that $\underline{P}(A)=\underline{P}(A_*)$. Then, the explicit
expression for $\underline{P}(A)$ is
$ \underline{P}(A_*)=\sum_{ E \subseteq
A} \underline{P}(E),$
which remains quite simple to compute, and more efficient than
computing the belief degree by checking which focal elements are
included in $A$.

Notice that Equation
(\ref{eq:LowerprobGenPboxes}) can be restated in terms of the two
possibility distributions $\pi_{\overline{F}},\pi_{\underline{F}}$,
rewriting $\underline{P}(E)$ as
$$\underline{P}(E)=\max(0,N_{\pi_{\underline{F}}}(\bigcup_{k=1}^{j}
G_k) - \Pi_{\pi_{\overline{F}}}(\bigcup_{k=1}^{i-1} G_k)),$$ where
$N_{\pi_i}(A),\Pi_{\pi_i}(A)$ are respectively the necessity and
possibility degree of event $A$ (given by Equations
(\ref{eq:possnecdegree})) with respect to a distribution $\pi_i$. It
makes $\underline{P}(A_*)$ even easier to compute.


\subsection{Probability intervals and generalized p-boxes}
\label{sec:GenPboxProbInt}

As in the case of random sets, there is no direct relationship
between probability intervals and generalized p-boxes. The two
representations have comparable complexities, but do not involve the
same kind of events. Nevertheless, given previous results, we can
state how a probability interval $L$ can be transformed into
a generalized p-box $[\underline{F},\overline{F}]$, and vice-versa.

First consider a  probability interval $L$ and some indexing
$\{x_1,\ldots,x_n\}$ of elements in $X$.
A generalized p-box $[\underline{F}',\overline{F}']$
outer-approximating the probability interval $L$ can be
computed by means of Equations (\ref{eq:EventProbInt}) as follows:
\begin{align}
& \underline{F}'(x_i) = \underline{P}(A_i)=\alpha'_i=\max(\sum_{x_i
\in A_i} l(x_i), 1 -
\sum_{x_i \notin A_i} u(x_i) ) \label{eq:InttoPbox} \\
& \overline{F}'(x_i) = \overline{P}(A_i)=\beta'_i =\min(\sum_{x_i
\in A_i} u(x_i), 1 - \sum_{x_i \notin A_i} l(x_i)) \nonumber
\end{align}
where $\underline{P},\overline{P}$ are respectively the lower and
upper probabilities of $\P_L$. Recall that $A_i=\{x_1,\ldots,x_i\}$.
Note that each permutation of elements of $X$ provide a different
generalized p-box and that there is no tightest outer-approximation
among them.

Next, consider a generalized p-box $[\underline{F},\overline{F}]$
with nested sets $A_1 \subseteq \ldots \subseteq A_n$. The 
probability interval $L'$ on elements $x_i$ corresponding to
$[\underline{F},\overline{F}]$ is given by:
\begin{align}
& \underline{P}(\{x_i\}) = l'(x_i) = \max(0,\alpha_i - \beta_{i-1})  \label{eq:PboxtoInt}\\
& \overline{P}(\{x_i\}) = u'(x_i) = \beta_i - \alpha_{i-1}
\nonumber
\end{align}
where $\underline{P},\overline{P}$ are the lower and upper
probabilities of the credal set $\P_{[\underline{F},\overline{F}]}$
on elements of $X$, $\alpha_i=\overline{F}(A_i)$,
$\beta_i=\underline{F}(A_i)$ and $\beta_0=\alpha_0=0$. This is the
tightest probability interval outer-approximating the
generalized p-box, and there is only such set.

Of course, transforming a probability interval $L$ into a
p-box $[\underline{F},\overline{F}]$ and vice-versa generally
induces a loss of information. But we can show that probability
intervals are representable (see definition~\ref{def:compa}) by
generalized p-boxes: let $\Sigma_{\sigma}$ the set of all possible
permutations $\sigma$ of elements of $X$, each defining a linear
order. A generalized p-box according to permutation $\sigma$ is
denoted $[\underline{F}',\overline{F}']_{\sigma}$ and called a
$\sigma$-p-box. We then have the following proposition (whose proof is in the appendix):
\begin{prop}
\label{prop:PBoxIntRel} Let $L$ be a probability interval,
and let
$[\underline{F}',\overline{F}']_{\sigma}$
be the  $\sigma$-p-box obtained from $L$ by applying Equations
(\ref{eq:InttoPbox}). Moreover, let
$L''_{\sigma}$ denote the probability interval obtained from
the $\sigma$-p-box $[\underline{F}',\overline{F}']_{\sigma}$  by
applying Equations (\ref{eq:PboxtoInt}). Then, the various credal
sets thus defined satisfy the following property:
\begin{equation}
\P_L = \bigcap_{\sigma \in \Sigma_{\sigma}}
\P_{[\underline{F}',\overline{F}']_{\sigma}} = \bigcap_{\sigma \in
\Sigma_{\sigma}} \P_{L''_{\sigma}}
\end{equation}
\end{prop}
This means that the information modeled by a set $L$ of probability
intervals can be entirely recovered by considering sets of
$\sigma$-p-boxes. Note that not all $|X|!$ such permutations need to
be considered, and that in practice, $L$ can be exactly recovered by means of
a reduced set $\mathcal{S}$ of $|X|/2$ permutations, provided that
$\{x_{\sigma(1)}, \sigma \in \mathcal{S} \}\cup \{x_{\sigma(n)},
\sigma \in \mathcal{S} \} = X$.
Since $\P_{[\underline{F},\overline{F}]}=\P_{\pi_{\underline{F}}}
\cap \P_{\pi_{\overline{F}}}$, then  it is immediate from
Proposition~\ref{prop:PBoxIntRel},that, in terms of credal sets,
$\P_L = \bigcap_{\sigma \in \Sigma_{\sigma}} \left(
\P_{\pi_{\underline{F}_{\sigma}}} \cap
\P_{\pi_{\overline{F}_{\sigma}}} \right) $, where
$\pi_{\underline{F}_{\sigma}},\pi_{\overline{F}_{\sigma}}$ are
respectively the possibility distributions corresponding to
$\underline{F}_{\sigma}$ and $\overline{F}_{\sigma}$.

\begin{figure}
\begin{center}
\begin{tikzpicture}
 \node (ImpProb) at (0,-0.8) [draw] {\scriptsize{Lower/upper prev.}};
 \node (LowUpp) at (0,-1.8) [draw] {\scriptsize{Lower/upper prob.}};
 \node (2mon) at (0,-2.8) [draw] {\scriptsize{2-monotone capacities}} ;
 \node (randset) at (0,-3.8) [draw] {\scriptsize{Random sets ($\infty$-monot)}} ;
 \node (GenPbox) at (0,-5.6) [draw,ultra thick] {\scriptsize{Generalized p-boxes}} ;
 \node (Pbox) at (0,-6.8) [draw] {\scriptsize{P-boxes}} ;
 \node (Prob) at (0,-7.8) [draw] {\scriptsize{Probabilities}} ;
 \node (ProbInt) at (-3,-4.6) [draw] {\scriptsize{Probability
 Intervals}} ;
 \node (Poss) at (3,-7.6) [draw] {\scriptsize{Possibilities}} ;

 \draw[->] (ImpProb) -- (LowUpp)  ;
 \draw[->]  (LowUpp) -- (2mon) ;
 \draw[->]  (2mon) -- (randset)  ;
 \draw[->,ultra thick]  (randset) -- (GenPbox) ;
 \draw[->,ultra thick] (GenPbox) -- (Pbox)  ;
 \draw[->]  (Pbox) -- (Prob) ;
 \draw[->,ultra thick] (GenPbox.south east) -- (Poss) ;
 \draw[->] (2mon.south west) -- (ProbInt)  ;
 \draw[->]  (ProbInt) -- (Prob.north west) ;
 \draw[dashed,->,ultra thick]  (Poss) .. controls
+(-0.2,1.2) .. (GenPbox.east) ; \draw[dashed,->,ultra thick]
(GenPbox) .. controls +(-0.4,1.2) .. (ProbInt) ;

\end{tikzpicture}  \caption{Representation relationships: summary with generalized p-boxes. $A \longrightarrow B$: A generalizes B. $A \dashrightarrow B$: B is representable by A }\label{fig:summarypre}
\end{center}
\end{figure}
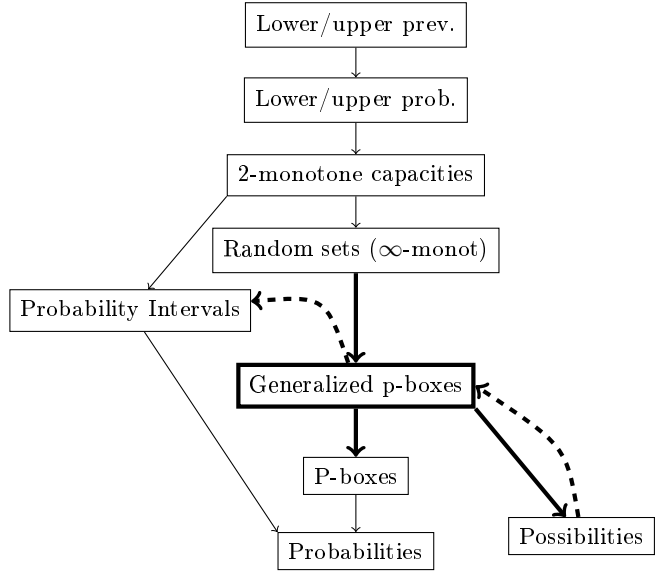

\section{Conclusion}



This paper introduces a generalized notion of p-box. Such a
generalization allows to define p-boxes  on finite
(pre)-ordered spaces as well as discretized p-boxes on
multi-dimensional spaces equipped with an arbitrary (pre)-order. On the real line, this preorder can be of the form $x \leq_\rho y$ if and only if $\vert x -\rho \vert \leq \vert y -\rho \vert$, thus  accounting for events of the form ``close to a prescribed  value $\rho$''.
Generalized
p-boxes are representable by a pair of comonotone possibility
distributions. They are special case of random sets, and the corresponding mass assignment has been laid bare. Generalized
p-boxes are thus more expressive than single possibility
distributions and likely to be more tractable than general random
sets. Moreover, the fact that they can be interpreted as lower and
upper confidence bounds over nested sets makes them quite attractive
tools for subjective elicitation. Finally, we showed the relation
existing between generalized p-boxes and sets of probability
intervals. Figure~\ref{fig:summarypre} summarizes the results of
this paper, by placing generalized p-boxes inside the graph of
Figure~\ref{fig:summarystate}. New relationships and representations
obtained in this paper are in bold lines.
Computational aspects of calculations with generalized p-boxes need
to be explored in greater detail (as is done by De Campos \emph{et
al.}~\cite{CamposAll94} for probability intervals) as well as their application to the elicitation of imprecise probabilities.
 Another
issue is to extend presented results to more general spaces, to
general lower/upper previsions or to cases not considered here (e.g.
continuous p-boxes with discontinuity points), possibly using
existing results~\cite{Smets05,GDecoo05}.
Interestingly, the key condition when representing generalized
p-boxes by two possibility distributions is their comonotonicity. In
the second part of this paper, we pursue the present study by dropping this
assumption. We then recover so-called  clouds, recently
proposed by Neumaier~\cite{Neumaier04}.


\section*{Appendix}

\begin{proof}[\textbf{Proof of Proposition~\ref{prop:PBoxBel}}]

From the nested sets $A_1 \subseteq A_2 \subseteq \ldots \subseteq
A_n=X$ we can build a partition s.t. $G_1=A_1, G_2=A_2\setminus A_1,
\ldots, G_n=A_n\setminus A_{n-1}$. Once we have a finite partition,
every possible set $B \subseteq X$ can be approximated from above
and from below by pairs of sets $ B_* \subseteq
B^*$~\cite{Pawlak91}:
 \begin{displaymath}
B^* = \bigcup \{G_i, G_i \cap B \neq \emptyset\};  \quad  B_* = \bigcup \{G_i, G_i \subseteq B\}
\end{displaymath}
made of a finite union of the partition elements intersecting  or
contained in this set $B$. Then $\underline{P}(B) =
\underline{P}(B_*)$,$\overline{P}(B)=\overline{P}(B^*)$, so we only
have to care about unions of elements $G_i$ in the sequel.
Especially, for each event $B \subset G_i$ for some $i$, it is clear
that $\underline{P}(B)=0=Bel(B)$ and
$\overline{P}(B)=\overline{P}(G_i)=Pl(B)$. So, to prove
Proposition~\ref{prop:PBoxBel}, we have to show that lower
probabilities given by a generalized p-box
$[\underline{F},\overline{F}]$ and by the corresponding random set
built through algorithm~\ref{alg:PboxBel} coincide on unions of
elements $G_i$. We will first concentrate on unions of conscutive
elements $G_i$, and then to any union of such elements.

Let us first consider union of consecutive elements
$\bigcup_{k=i}^{j} G_k$ (when $k=1$, we retrieve the sets $A_j$).
Finding $\underline{P}(\bigcup_{k=i}^{j} G_k)$ is equivalent to
computing the minimum of $ \sum_{k=i}^{j} P(G_k)$ under the
constraints
\begin{displaymath}
i=1,\ldots,n \qquad \alpha_i \leq P(A_i)=\sum_{k=1}^{i} P(G_k) \leq
\beta_i
\end{displaymath}
which reads
\begin{displaymath}
\alpha_j \leq P(A_{i-1}) + \sum_{k=i}^{j} P(G_k) \leq \beta_j
\end{displaymath}
so $\sum_{k=i}^{j} P(G_k) \geq \max(0,\alpha_j - \beta_{i-1})$. This
lower bound is optimal, since it is always reachable: if $\alpha_j
> \beta_{i-1}$, take $P$ s.t. \mbox{$P(A_{i-1}) =
\beta_{i-1}$}, \mbox{$P(\bigcup_{k=i}^{j} G_k) = \alpha_j -
\beta_{i-1}$}, \mbox{$P(\bigcup_{k=j+1}^{n} G_k) =  1- \alpha_j.$}
If $\alpha_j \leq \beta_{i-1}$, take $P$ s.t. \mbox{$P(A_{i-1}) =
\beta_{i-1}$,} \mbox{$P(\bigcup_{k=i}^{j} G_k) = 0$},
\mbox{$P(\bigcup_{k=j+1}^{n} E_k) =  1- \beta_{i-1}.$} And we can see, by looking at Algorithm~\ref{alg:PboxBel}, that
$Bel(\bigcup_{k=i}^{j} G_k)=\max(0,\alpha_j - \beta_{i-1})$: focal
elements of $Bel$ are subsets of $\bigcup_{k=i}^{j} G_k$ if
$\beta_{i-1} < \alpha_j$ only.


Now, let us consider a union $A$ of non-consecutive elements s.t. \\
\mbox{$A = (\bigcup_{k=i}^{i+l} G_k \ \cup \ \bigcup_{k=i+l+m}^{j}
G_k)$} with $m > 1$. As in the previous case, we must compute
\mbox{$\min \Big( \sum_{k=i}^{i+l} P(G_k) + \sum_{k=i+l+m}^{j}
P(G_k) \Big)$} to find the lower probability on $\underline{P}(A)$.
An obvious lower bound is given by
\begin{displaymath}
\min \Big( \sum_{k=i}^{i+l} P(G_k) + \sum_{k=i+l+m}^{j} P(G_k) \Big)
\geq \min \Big( \sum_{k=i}^{i+l} P(G_k) \Big) + \min \Big(
\sum_{k=i+l+m}^{j} P(G_k) \Big)
\end{displaymath}
and this lower bound is equal to
\begin{displaymath}
\max(0,\alpha_{i+l} - \beta_{i-1}) + \max(0,\alpha_j -
\beta_{i+l+m-1}) = Bel(A)
\end{displaymath}
Consider the two following cases with probabilistic mass assignments
showing that bounds are attained:
\begin{itemize}
\item $\alpha_{i+l} < \beta_{i-1}$, $\alpha_j < \beta_{i+l+m-1}$ and probability masses:\\  \mbox{$P(A_{i-1}) = \beta_{i-1}$}, \mbox{$P(\bigcup_{k=i}^{i+l} G_k) = \alpha_{i+l} - \beta_{i-1}$}, \mbox{$P(\bigcup_{k=i+l+1}^{i+l+m-1} G_k) = \beta_{i+l+m-1} - \alpha_{i+l}$}, \\ \mbox{$P(\bigcup_{k=i+l+m}^{j} G_k) = \alpha_j - \beta_{i+l+m-1}$} and \mbox{$P(\bigcup_{k=j+1}^{n} G_k) = 1 -
\alpha_j$.}
\item $\alpha_{i+l} > \beta_{i-1}$, $\alpha_j > \beta_{i+l+m-1}$ and probability masses:\\  \mbox{$P(A_{i-1}) = \beta_{i-1}$}, \mbox{$P(\bigcup_{k=i}^{i+l} G_k) = 0$}, \mbox{$P(\bigcup_{k=i+l+1}^{i+l+m-1} G_k) = \alpha_j - \beta_{i-1}$},\\  \mbox{$P(\bigcup_{k=i+l+m}^{j} E_k) = 0$} and \mbox{$P(\bigcup_{k=j+1}^{n} G_k) = 1 -
\alpha_j$.}
\end{itemize}
The same line of thought can be followed for the two remaining cases.
As in the consecutive case, the lower bound is reachable without
violating any of the restrictions associated to the generalized
p-box. We have $\underline{P}(A)= Bel(A)$ and the extension of this
result to any number $n$ of "discontinuities" in the sequence of
$G_k$ is straightforward. The proof is complete, since for every possible union $A$ of
elements $G_k$, we have $\underline{P}(A) = Bel(A)$
\end{proof}

\begin{proof}[\textbf{Proof of Proposition~\ref{prop:PBoxIntRel}}]
To prove this proposition, we must first recall a result
 given by De Campos \emph{et al.}~\cite{CamposAll94}: given two probability intervals $L$ and
$L'$ defined on a space $X$ and the induced credal sets $\P_L$ and
$\P_{L'}$, the conjunction $\P_{L \cap L'}=\P_L \cap \P_{L'}$ of
these two sets can be modeled by the set $(L\cap L')$ of probability
intervals that is such that for every element $x$ of $X$,
$$l_{(L\cap L')}(x)=\max(l_{L}(x),l_{L'}(x)) \textrm{ and }
u_{(L\cap L')}(x)=\min(u_{L}(x),u_{L'}(x))$$ and these formulas
extend directly to the conjunction of any number of
probability intervals on $X$.

To prove Proposition~\ref{prop:PBoxIntRel}, we will show, by using
the above conjunction, that $\P_L = \bigcap_{\sigma \in
\Sigma_{\sigma}} \P_{L''_{\sigma}}$. Note that we have, for any
$\sigma \in \Sigma_{\sigma}$, $\P_L \subset
\P_{[\underline{F}',\overline{F}']_{\sigma}} \subset
\P_{L''_{\sigma}}$, thus showing this equality is sufficient to
prove the whole proposition.

Let us note that the above inclusion relationships alone ensure us
that \\ \mbox{$\P_L \subseteq \bigcap_{\sigma \in \Sigma_{\sigma}}
\P_{[\underline{F}',\overline{F}']_{\sigma}} \subseteq
\bigcap_{\sigma \in \Sigma_{\sigma}} \P_{L''_{\sigma}}$}. So, all we
have to show is that the inclusion relationship is in fact an
equality.

Since we know that both $\P_L$ and $\bigcap_{\sigma \in
\Sigma_{\sigma}} \P_{L''_{\sigma}}$ can be modeled by 
probability intervals, we will show that the lower bounds $l$ on
every element $x$ in these two sets coincide (and the proof for
upper bounds is similar).

For all $x$ in $X$, $l_{L''_{\Sigma}}(x)=\max_{\sigma \in
\Sigma_{\sigma}}\{l_{L''_{\sigma}}(x)\}$, with $L''_{\Sigma}$ the
probability interval corresponding to $\bigcap_{\sigma \in
\Sigma_{\sigma}} \P_{L''_{\sigma}}$ and $L''_{\sigma}$ the 
probability interval corresponding to a particular permutation
$\sigma$. We must now show that, for all $x$ in $X$,
$l_{L''_{\Sigma}}(x)=l_L(x)$.

Given a ranking of elements of $X$, and by applying successively
Equations (\ref{eq:InttoPbox}) and (\ref{eq:PboxtoInt}), we can
express the differences between bounds $l''(x_i)$ of the set $L''$
and $l(x_i)$ of the set $L$ in terms of set of bounds $L$. This
gives, for any $x_i \in X$:
\begin{align}\label{eq:IntDiff}  l(x_i) -
l''(x_i)=&\min(l(x_i),0 +\sum\limits_{x_i \in A_{i-1}} (u(x_i) -
l(x_i)),\\ &  0+ \sum\limits_{x_i \in A_{i}^c} (u(x_i) - l(x_i)) ,
(l(x_i) + \sum_{\substack{x_j \neq x_i \\ x_j \in X}} u(x_j))-1,
1-\sum_{x_i \in X} l(x_i)  )\nonumber \end{align}

We already know that, for any permutation $\sigma$ and for all $x$
in $X$, we have $l_{L''_{\sigma}}(x) \leq l_L(x)$. So we must now
show that, for a given $x$ in $X$, there is one permutation $\sigma$
such that $l_{L''_{\sigma}}(x) = l_L(x)$. Let us consider the
permutation placing the given element at the front. If $x$ is the
first element $x_{\sigma(1)}$, then Equation (\ref{eq:IntDiff}) has
value $0$ for this element, and we thus have $l_{L''_{\sigma}}(x) =
l_L(x)$. Since if we consider every possible ranking, every element
$x$ of $X$ will be first in at least one of these rankings, this
completes the proof.
\end{proof}

\end{document}